\documentclass{amsart}
\usepackage{amsmath}

\newcommand{\R}{{\mathbb R}}
\newcommand{\RP}{{\mathbb R}P}

\newcommand{\sys}{{\rm sys}}
\newcommand{\vol}{{\rm vol}}

\newcommand{\dmv}{\tilde{V}}  
\newcommand{\emv}{\tilde{W}}  
\newcommand{\rs}{\ \tilde{+}\ }   
\newcommand{\ctg}[1]{T^{*}#1} 
\newcommand{\csph}[1]{S^{*}#1}
\newcommand{\cod}[1]{D^{*}#1}

\newtheorem{theorem}{Theorem}[section]
\newtheorem{corollary}{Corollary}[section]
\newtheorem{lemma}{Lemma}[section]
\newtheorem{proposition}{Proposition}[section]
\newtheorem*{conjecture}{Uniformization conjecture}

\theoremstyle{definition}
\newtheorem{definition}{Definition}[section]

\title{Dual Mixed Volumes and Isosystolic Inequalities}

\author{J.C. \'Alvarez Paiva}

\address{J.C. \'Alvarez Paiva, Department of Applied Mathematics,
Polytechnic University, Six MetroTech Center, Brooklyn, NY, 11201,
USA.}

\email{jalvarez@duke.poly.edu}



\begin{document}%

\begin{abstract}
The theory of dual mixed volumes is extended to star bodies in
cotangent bundles and is used to prove several isosystolic inequalities
for Hamiltonian systems and Finsler metrics.
\end{abstract}
\maketitle

\tableofcontents

\section{Introduction}

A celebrated theorem of Pu (see \cite{Pu}) states that the volume of
any Riemannian metric on the projective plane is bounded below by
$2/\pi$ times the square of the length of the shortest non-contractible
geodesic. Equality holds if and only if the metric is of constant curvature.
In the same paper, Pu investigates analogues of this result in other homogeneous
spaces under the condition that the metrics involved be conformal to an invariant
metric. In \cite{Berger}, M. Berger considered infinitesimal deformations of
metrics in compact symmetric spaces of rank one and proved, among other results,
that if $g_t$ is a one-parameter family of Riemannian metrics on $\RP^n$ such
that $g_0$ is the standard invariant metric of constant curvature, then there
exists a second family of Riemannian metrics, $h_t$, that agrees to first order
with $g_t$ at $t = 0$ and satisfies the isosystolic inequality
\begin{equation}\label{Berger's_inequality}
\frac{\sys_{1}^{n}(\RP^n,h_t)}{\vol(\RP^n,h_t)} \leq
\frac{\sys_{1}^{n}(\RP^n,g_{0})}{\vol(\RP^n,g_{0})} \ .
\end{equation}
Here $\sys_{1}(\RP^n,g)$, the {\it $1$-systole\/} of $(\RP^n,g)$, denotes the
length of the shortest non-contractible geodesic for the metric $g$.

This paper extends many of Pu's sharp isosystolic inequalities and Berger's
infinitesimal isosystolic inequality~(\ref{Berger's_inequality}) to Finsler
metrics and Hamiltonian systems. In mechanical terms these inequalities provide
lower bounds for the Liouville volume enclosed by an energy surface in terms of
the action of periodic solutions of the system on that energy level.

The first part of the paper extends the theory of dual mixed volumes to star
bodies in cotangent bundles. This theory is used to define relative
invariants for pairs of Finsler metrics on a compact $n$-dimensional
manifold $M$, $\emv_k(M;L,L_0)$, $1 \leq k \leq n-1$, and to prove the
inequality
\begin{equation}\label{volume_inequality}
\emv_k(M;L,L_0)^n \leq  \frac{\vol(M,L)^{n-k}}{\vol(M,L_{0})^{n-k}} \ .
\end{equation}
Equality holds if and only if $L$ is a constant multiple of $L_0$.

The second part of the paper considers three different sets of hypotheses under
which the inequality
\begin{equation}\label{systole_inequality}
\emv_k(M;L,L_0) \geq \frac{\sys_{k}(M,L)}{\sys_{k}(M,L_{0})}
\end{equation}
holds. To recall, for $k > 1$, the {\it $k$-systole\/} of a Riemannian or
Finsler manifold $(M,L)$, $\sys_k (M,L)$, is the infimum of the areas of all
$k$-dimensional submanifolds that are not homologous to zero. The notion of
area and volume on a Finsler manifold we shall use is that of
Holmes and Thompson \cite{Holmes-Thompson, Alvarez-Thompson}, which from the
Hamiltonian viewpoint is more natural than the Hausdorff measure.

Under these hypotheses, inequalities~(\ref{volume_inequality})
and~(\ref{systole_inequality}) yield the isosystolic inequality
\begin{equation}\label{isosystolic_inequality}
\frac{\sys_{k}(M, L)^n}{\sys_{k}(M,L_{0})^n} \leq \emv_k(M;L,L_0)^n
\leq \frac{\vol(M,L)^{n-k}}{\vol(M,L_{0})^{n-k}} \ .
\end{equation}

In the first set of hypotheses $M$ is a homogeneous space, $L_0$ is an
invariant Finsler metric, $L$ is conformal to $L_0$, and the result is a
fairly straight forward generalization of the results in \cite{Pu}. The second
and third sets of hypotheses are more subtle.

\begin{theorem}\label{main:1}
Let $L_0$ be a Finsler metric on $\RP^n$ such that its geodesic flow is
symplectically conjugate to the geodesic flow of a metric of constant
curvature. If $L$ is conformal to $L_0$, then
$$
\frac{\sys_{1}(M, L)^n}{\sys_{1}(M,L_{0})^n} \leq \emv_{n-1}(M;L,L_0)^n
\leq \frac{\vol(M,L)}{\vol(M,L_{0})} \ .
$$
\end{theorem}

\begin{theorem}\label{main:2}
Let $L_0$ be a Finsler metric on $\RP^n$ such that its geodesic flow is periodic.
If the geodesic flow of a Finsler metric $L$ on $\RP^n$ commutes with that of
$L_0$, then
$$
\frac{\sys_{1}(M, L)^n}{\sys_{1}(M,L_{0})^n} \leq \emv_{n-1}(M;L,L_0)^n
\leq \frac{\vol(M,L)}{\vol(M,L_{0})} \ .
$$
\end{theorem}

It is not clear if there is any difference between the hypotheses
on $L_0$ in the previous theorems. The geodesic flow of a Finsler metric on
$\RP^2$ is periodic if and only if it is symplectically conjugate to the geodesic
flow of a metric of constant curvature. In higher dimensions this is also the case
for all known examples.

Theorem~\ref{main:2} together with the technique of averaging of Hamiltonian
systems (see \cite{Moser} and \cite{Cushman}) yields the following generalization
of Berger's result.

\begin{theorem}\label{main:3}
If $L_{t}$ is a smooth path of Finsler metrics on $\RP^n$ such
that the geodesic flow of $L_{0}$ is periodic, then there exists another
smooth path of Finsler metrics, $K_t$, that agrees to first order with $L_t$
at $t = 0$ and satisfies the isosystolic inequality
$$
\frac{\sys_{1}^{n}(\RP^n,K_t)}{\vol(\RP^n,K_t)} \leq
\frac{\sys_{1}^{n}(\RP^n,L_{0})}{\vol(\RP^n,L_{0})} \ .
$$
\end{theorem}

Unfortunately, the Finsler extension of Pu's theorem for Riemannian metrics on
the projective plane lies beyond the reach of these techniques. This extension
is, however, an easy consequence of the main result of S. Ivanov in \cite{Ivanov}
and would follow immediately from Theorem~\ref{main:1} if the following Finsler
generalization of the uniformization theorem were true:

\begin{conjecture}
If $L$ is a reversible Finsler metric on $\RP^2$, then there
exists a smooth positive function $\rho$ such that the geodesic
flow of the metric $\rho L$ is periodic.
\end{conjecture}

It is important to remark that the extensive work on coarse isosystolic
inequalities and systolic freedom by Babenko, Gromov, Katz, and Suciu
(see, for example, \cite{Gromov:1,Gromov:2,Babenko-Katz,Katz-Suciu})
applies unchanged to Finsler metrics.

\medskip
\noindent {\bf Acknowledgements.}
The author warmly thanks Erwin Lutwak for having introduced him to dual mixed
volumes and is grateful to Misha Katz and Deane Yang for their comments on an
earlier version of this work.

\section{Dual mixed volumes on cotangent bundles}

The theory of dual mixed volumes was introduced by E. Lutwak in
\cite{Lutwak:1, Lutwak:2} as a version of the Brunn-Minkowski
theory in which averages of areas of projections of convex bodies,
{\it quermassintegrals,} are replaced by averages of areas of
central sections of star-shaped bodies. The theory plays a
key role in the solution of the Busemann-Petty problem and
other problems in geometric tomography (\cite{Lutwak:3,Gardner}).
It also has applications to integral geometry (e.g., \cite{Zhang})
and the theory of valuations (\cite{Klain:1,Klain:2}).

In this section we present an straight-forward extension of the dual theory
to star bodies in cotangent bundles. In simple terms, a star body in a cotangent
bundle is a choice of a star-shaped body in each cotangent space that varies
continuously with the base point. The following is the precise definition:

\begin{definition}
Let $M$ be a compact manifold whose boundary need not be empty. A
{\it star Hamiltonian\/} is a continuous function
$H : \ctg{M} \rightarrow [0,\infty)$ that is
\begin{itemize}
\item positive outside the zero section;
\item positively homogeneous of degree one (i.e., $H(tp) = tH(p)$
      whenever $t > 0$);
\item proper.
\end{itemize}
We shall say a star Hamiltonian is {\it smooth\/} if it is smooth outside the
zero section.
A subset $A \subset \ctg{M}$ is said to be a {\it star body \/} if
$$
A = \cod{(M,H)} := \{p \in \ctg{M} : H(p) \leq 1 \}
$$
for some star Hamiltonian.
\end{definition}

Examples of star bodies are the unit codisc bundles of Riemannian or
Finsler metrics and, at least in this paper, the most interesting applications
will deal with these.

A very useful tool in the study of star bodies is the {\it radial function:}
if $A$ is a star body with Hamiltonian $H$ and $U$ is a ``model''
star body with Hamiltonian $H_{0}$ and boundary $\partial U$, the radial
function, $\rho_{A}$, of $A$ (with respect to $U$) is the
restriction of $1/H$ to $\partial U$. It's easy to see that any continuous,
positive function on $\partial U$ is the radial function of some star body,
and that the map   $p \mapsto \rho_{A}(p)\, p$
is a homeomorphism between the boundary of $U$ and the boundary of $A$.

In the classical theory of dual mixed volumes, the model star-shaped body is the
unit sphere. However, most of the basic results are independent of the choice
of model body, and in great part the flexibility of the extension of the
theory to cotangent bundles is due to this. Most of the results in this paper
are obtained by comparing star bodies to different, more symmetric bodies.

An elementary property of star bodies is that they form a lattice: finite unions
and intersections of star bodies are star bodies. Indeed, the radial function
for the union of two star bodies is the maximum of their radial functions, and
the radial function of the intersection of two star bodies is the minimum of
their radial functions. More importantly, we can dilate and add star bodies:

\begin{definition}
Let $A$ and $B$ be star bodies with radial functions $\rho_{A}$ and $\rho_{B}$.
If $\lambda$ is a positive real number, we denote by $\lambda A$ the star body
with radial function $\lambda \rho_{A}$ and by $A \rs B$, the {\it radial sum\/}
of $A$ and $B$, the star body with radial function $\rho_{A} + \rho_{B}$.
\end{definition}

An intrinsic, more geometric description of the radial sum goes as follows:
if $p$ and $p^{\prime}$ are two covectors in $T^{*}_{x}{M}$, we define
their {\it radial sum,} $p \rs p^{\prime}$,
as $p + p^{\prime}$ if they belong to the same one-dimensional subspace
and as zero otherwise. The radial sum, $A \rs B$, of two star bodies $A$ and $B$ in
$\ctg{M}$ is the union over all points $x$ in $M$ of the sets
$$
(A \rs B)_{x} := \{p \rs p^{\prime} :
 p \in A \cap T_{x}^{*}M, \ p^{\prime} \in B \cap T_{x}^{*}M \}.
$$

Radial functions are also useful in describing a topology on the set of star
bodies.

\begin{definition}
Let us fix a model body and define the distance between two star bodies
$A$ and $B$ with radial functions $\rho_{A}$ and $\rho_{B}$ as the maximum of
$|\rho_{A} - \rho_{B}|$. The topology induced by this metric is independent of
the choice of model body used to define the radial functions, and will be called
the {\it radial Hausdorff topology.}
\end{definition}

The theory of dual mixed volumes is the study of the interaction between the
radial sum and the volume of star bodies. An advantage of extending the classical
theory to cotangent bundles is that we may use the natural symplectic volume
on these spaces:

Let $M$ be a smooth $n$-dimensional manifold and let $\pi : \ctg{M} \rightarrow M$
be its  cotangent bundle. The {\it canonical 1-form\/} $\alpha$ on $\ctg{M}$
is the form whose value at a tangent vector $v_{p} \in T_{p}\ctg{M}$
equals $p(\pi_{*} (v_{p}))$. In local canonical coordinates
$(q_{1},\dots,q_{n},p_{1},\dots,p_{n})$, $\alpha$ takes the form
$$
\alpha = \sum_{i=1}^{n} p_{i}dq_{i} \ .
$$
The {\it symplectic form\/} on $\ctg{M}$ is defined as the 2-form
$\omega := d \alpha$ and the (Liouville) {\it volume form\/} is
$\omega^{n}/n!$.

Notice that when $(M,g)$ is an $n$-dimensional Riemannian manifold
the volume of its unit codisc bundle equals the Riemannian volume of  $(M,g)$
times the volume of the Euclidean unit ball of dimension $n$.

The volume of a star body can be easily described in terms of its radial function:

\begin{proposition}
Let $A$ be a star body in the cotangent bundle of an $n$-dimensional manifold
$M$ and let $\rho_{A}$ be its radial function. If we set
$\Omega := \alpha \wedge (d\alpha)^{n-1}/n!$, the volume of $A$, $V(A)$, is given
by the integral
$$
\int_{\partial U} \! \rho_{A}^{n} \ \Omega .
$$
\end{proposition}

\proof
We give a proof in the case where both the model body $U$ and the star body $A$
have smooth boundaries. The general result follows by a standard approximation
argument.

By Stokes formula we have that
$$
V(A) := \int_{A} \omega^{n}/n! = \int_{\partial A} \Omega .
$$
The map $\delta : \partial U \rightarrow \partial A$ defined by
$\delta(p) = \rho(p)p$ is a diffeomorphism and since
$\delta^{*}\alpha = \rho \alpha$, we have that $\delta^{*}\Omega$
equals $\rho^{n} \Omega$. We conclude
that
\begin{eqnarray*}
V(A) & = & \int_{\partial A} \Omega = \int_{\delta(\partial U)} \Omega \\
     & = & \int_{\partial U} \delta^* \Omega =  \int_{\partial U} \! \rho^{n} \Omega .
\end{eqnarray*}
\qed

\begin{definition}
Let $A_{1},\dots,A_{n}$ be $n$ star bodies in the cotangent of a compact
$n$-dimensional manifold $M$ and let $\rho_{1},\dots,\rho_{n}$ be their radial
functions. The {\it dual mixed volume\/} of $A_{1},\dots,A_{n}$, denoted by
$\dmv(A_{1},\dots,A_{n})$, is defined as the integral
$$
\dmv(A_{1},\dots,A_{n}) := \int_{\partial U} \! \rho_{1} \cdots \rho_{n} \ \Omega \ .
$$
\end{definition}

The next proposition shows that dual mixed volumes do not depend on the choice
of model body that is used in the definition of the radial functions.

\begin{proposition}
Let $A_{1},\dots,A_{k}$ be star bodies in the cotangent bundle of an
$n$-dimensional compact manifold $M$. The volume of
$\lambda_{1}A_{1} \rs \cdots \rs \lambda_{k}A_{k}$ is an
$n$-th-degree polynomial in the $\lambda_{i}$,
$$
V(\lambda_{1}A_{1} \rs \cdots \rs \lambda_{k}A_{k}) =
\sum \lambda_{i_1}\cdots\lambda_{i_n} \dmv(A_{i_1},\dots,A_{i_n}) ,
$$
where the sum is taken over all positive integers less than $k$.
\end{proposition}

\proof
\begin{eqnarray*}
V(\lambda_{1}A_{1} \rs \cdots \rs \lambda_{k}A_{k}) & = &
\int_{\partial U} \!
\left(\lambda_{1}\rho_{1} + \cdots + \lambda_{k}\rho_{k} \right)^{n} \ \Omega  \\
& = &
\int_{\partial U} \left( \sum \lambda_{i_1} \cdots \lambda_{i_n}
                     \rho_{i_1} \cdots \rho_{i_n} \right) \ \Omega \\
& = &
\sum \lambda_{i_1}\cdots\lambda_{i_n} \dmv(A_{i_1},\dots,A_{i_n}) \ .
\end{eqnarray*}
\qed

Some of the basic properties of dual mixed volumes are:
\begin{itemize}
\item Continuity with respect to the radial Hausdorff topology.
\item Positivity: $\dmv(A_{1},\dots,A_{n}) > 0$.
\item Homogeneity:
$\dmv(\lambda_{1}A_{1},\dots,\lambda_{n}A_{n}) = \lambda_{1}\cdots\lambda_{n}
\dmv(A_{1},\dots,A_{n})$, $\lambda_{i} > 0$.
\item Strict monotonicity: if $A_{i} \subset B_{i}$ for all $i$, then
$$
\dmv(A_{1},\dots,A_{n}) \leq \dmv(B_{1},\dots,B_{n}) \ .
$$
Equality holds if and only if $A_{i}=B_{i}$ for all $i$.
\item $\dmv(A,\dots,A) = V(A)$.
\end{itemize}

Some useful abbreviations are
$$
\dmv_{k}(A,B) := \dmv(\underbrace{A,\dots,A}_{n-k} , \underbrace{B,\dots,B}_{k})
\quad \mbox{and} \quad \emv_{k}(A) := \dmv_{k}(A,U)/V(U) \ .
$$
In particular, notice that $\emv_{k}(A)$ is just the average of the
$(n-k)$-th power of the radial function of $A$:
$$
\emv_{k}(A) = \frac{1}{V(U)} \int_{\partial U} \! \rho_{A}^{n-k} \ \Omega \ .
$$

One of the most basic results about dual mixed volumes is the following
inequality:

\begin{theorem}[Main inequality]
If $A_{1},\dots,A_{n}$ are star bodies in the cotangent bundle of an
$n$-dimensional compact manifold,
$$
\dmv(A_{1},\dots,A_{n})^n \leq V(A_{1})\cdots V(A_{n}) \ .
$$
Moreover, the equality if and only if all the star bodies are dilations of
each other.
\end{theorem}

\proof
In terms of the radial functions $\rho_{1},\dots,\rho_{n}$ of the star bodies
$A_{1},\dots A_{n}$, we must show that
$$
\left( \int_{\partial U} \! \rho_{1}\cdots\rho_{n} \ \Omega \right)^{n} \
\leq \int_{\partial U}\! \rho_{1}^{n} \ \Omega \
\cdots \int_{\partial U}\! \rho_{n}^{n} \ \Omega \ .
$$
To do this set $\nu_{i} := \rho_{i}^{n}$, $1 \leq i \leq n$, and
consider the quantity
$$
\frac{\int \! \sqrt[n]{\nu_{1}\cdots\nu_{n}} \ \Omega}{\left(
\int \nu_{1} \ \Omega \right)^{1/n} \cdots  \
\left(\int \nu_{n} \ \Omega \right)^{1/n}}
=
\int \left(\frac{\nu_{1}}{\int \nu_{1} \ \Omega}\right)^{1/n} \cdots
      \left(\frac{\nu_{n}}{\int \nu_{n} \ \Omega}\right)^{1/n} \Omega \ ,
$$
where we have suppressed the region of integration, $\partial U$, to simplify the
notation. By the arithmetic-geometric mean inequality, this quantity is less than
$$
\int \frac{1}{n} \left( \frac{\nu_{1}}{\int \nu_{1} \ \Omega} + \cdots +
\frac{\nu_{n}}{\int \nu_{n} \ \Omega} \right) \ \Omega = 1 .
$$
Equality holds if and only if all radial functions are multiples of each other.
The result follows.
\qed

\begin{corollary}[Dual Minkowski inequalities]
If $A$ and $B$ are star bodies in the cotangent bundle of an
$n$-dimensional compact manifold,
$$
\dmv_{1}(A,B)^n \leq V(A)^{n-1} V(B) \quad \mbox{and} \quad
\dmv_{n-1}(A,B)^n \leq V(A)V(B)^{n-1} .
$$
\end{corollary}

\begin{theorem}[Dual Brunn-Minkowski inequality]
If $A$ and $B$ are star bodies in the cotangent bundle of an
$n$-dimensional compact manifold,
$$
V(A \rs B)^{1/n} \leq V(A)^{1/n} + V(B)^{1/n} \ .
$$
Equality holds if and only if the star bodies are dilations of each other.
\end{theorem}

\proof
Let $A_{1},A_{2}$, and $B$ be star bodies. Using the additivity of
dual mixed volumes and the first of the dual Minkowski inequalities,
we have that
\begin{eqnarray*}
\dmv_{1}(B,A_{1} \rs A_{2}) & = & \dmv_{1}(B,A_{1}) + \dmv_{1}(B,A_{2}) \\
& \leq & V(B)^{(n-1)/n}
         \left(V(A_{1})^{1/n} + V(A_{2})^{1/n}\right) .
\end{eqnarray*}
In the particular case when $B = A_{1} \rs A_{2}$, we also have that
left hand side of the previous inequality, $\dmv_{1}(B, A_{1} \rs A_{2})$,
is equal to $V(A_{1} \rs A_{2})$ and the dual Brunn-Minkowski inequality
follows immediately.
\qed

\section{Invariance of dual mixed volumes}

A basic remark in the classical theory is that dual mixed volumes
are invariant under the special linear group. This is because
special linear transformations preserve both radial sums and volumes.
However, something that seems to have escaped notice until now is that the
symmetry group is actually much larger: any volume-preserving transformation
of $\R^n$ that is positively homogeneous of degree one preserves dual mixed
volumes. In this section we extend this remark to our generalized setting.

\begin{definition}
A diffeomorphism $\phi : \ctg{M} \setminus 0 \rightarrow \ctg{M} \setminus 0$
is said to be a {\it special homogeneous transformation\/} if it is positively
homogeneous of degree one (i.e., $\phi(tp) = t\phi(p)$, $t > 0$) and
preserves the form $\Omega$. Equivalently, a special homogeneous transformation
is a volume-preserving transformation that is positively homogeneous of degree
one.
\end{definition}

\begin{theorem}\label{rs-volume}
Let $M$ be a compact manifold and let
$\phi : \ctg{M} \setminus 0 \rightarrow \ctg{M} \setminus 0$
be a special homogeneous transformation.
If $A$ and $B$ are star bodies in the cotangent bundle of $M$, then
$$
V(A \rs B) = V(\phi(A) \rs \phi(B)) \ .
$$
\end{theorem}

The proof of this theorem depends on the following trivial lemma:

\begin{lemma}
Let $\phi : \ctg{M} \setminus 0 \rightarrow \ctg{M} \setminus 0$ be a
diffeomorphism that is positively homogeneous of order one.
If $\rho_{A}$ is the radial function of a star body $A$ with respect to the
model body $U$, then $\rho_{A} \circ \phi^{-1}$ is the radial function of
$\phi(A)$ with respect to $\phi(U)$.
\end{lemma}

\noindent
{\it Proof of theorem \ref{rs-volume}.\ }
Using the lemma and the formula for the volume of the radial sum
$\phi(A) \rs \phi(B)$ in terms of its radial function with respect to $\phi(U)$,
we have that
$$
V(\phi(A) \rs \phi(B)) =  \int_{\phi(\partial U)}
(\rho_{A} \circ \phi^{-1} + \rho_{B} \circ \phi^{-1})^n \ \Omega \ .
$$
Since $\phi^{-1*} \Omega = \Omega$, we may write this integral as
$$
\int_{\phi(\partial U} \phi^{-1*}\left[(\rho_{A} + \rho_{B})^{n} \Omega\right]
= \int_{\partial U}(\rho_{A} + \rho_{B})^{n} \ \Omega = V(A \rs B) \ .
$$
\qed

\begin{theorem}
Let $A_{1},\dots,A_{n}$ and $B_{1},\dots,B_{n}$ be star bodies in the cotangent
bundle of a compact $n$-dimensional manifold $M$. If there exists a special
homogeneous transformation
$\phi : \ctg{M} \setminus 0 \rightarrow \ctg{M} \setminus 0$
such that $\phi(A_{i}) \subset B_{i}$, for all $i$, then
$$
\dmv(A_{1},\dots,A_{n}) \leq \dmv(B_{1},\dots,B_{n}) \ .
$$
\end{theorem}

\proof
Because of the basic monotonicity property of dual mixed volumes, we only need to
show that
$$
\dmv(\phi(A_{1}),\dots,\phi(A_{n}) = \dmv(A_{1},\dots,A_{n}) \ .
$$
To do this, notice that for any positive numbers $\lambda_{1},\dots,\lambda_{n}$,
\begin{eqnarray*}
\sum \lambda_{i_1}\cdots \lambda_{i_n} \ \dmv(A_{i_1},\dots,A_{i_n}) & = &
V(\lambda_{1}A_{1} \rs \cdots \rs \lambda_{n}A_{n}) \\
& = & V(\phi(\lambda_{1}A_{1} \rs \cdots \rs \lambda_{n}A_{n})) \\
& = & V(\lambda_{1}\phi(A_{1}) \rs \cdots \rs \lambda_{n}\phi(A_{n})) \\
& = &
\sum \lambda_{i_1}\cdots \lambda_{i_n} \ \dmv(\phi(A_{i_1}),\dots,\phi(A_{i_n}))
\ .
\end{eqnarray*}
This immediately implies the invariance of dual mixed volumes under
special homogeneous transformations.
\qed

Let us remark that $\dmv(A_{1},\dots,A_{n})$ is a non trivial invariant of
$A_{1},\dots,A_{n}$ in the sense that it is not a function of their
symplectic volumes. In fact, assume that the $A_{i}$'s are not dilations of each
other and that all have unit volume. If $B$ is another star body with unit volume,
the main inequality tells us that
$$
\dmv(A_{1},\dots,A_{n}) < 1 = \dmv(B,\dots,B) \ ,
$$
and so the dual mixed volume cannot be a function of the symplectic volume.

An important class of special homogeneous transformations are diffeomorphisms
from $T^*M \setminus 0$ to itself that preserve the canonical form $\alpha$.
These transformations are known as {\it homogeneous canonical transformations.}

While the results in this section show that dual mixed volumes are invariant
under homogeneous canonical transformations, their invariance under the larger
special group of homogeneous transformations makes them uninteresting as
symplectic invariants.

\section{Finsler metrics and optical Hamiltonians}

In this section, we will quickly review the basic concepts in Finsler geometry
that will be used throughout the rest of the paper.

Roughly speaking, a Finsler metric on a manifold $M$ is a smooth choice of norm
on each tangent space of $M$. However, in this paper we shall sometimes work with
a slightly more general geometric structure where the norm of a vector $v$ is not
necessarily equal to the norm of $-v$. Moreover, the unit sphere of this
``non-symmetric" norm on each tangent space $T_x M$ must be
{\it quadratically convex:} its principal curvatures are positive for any auxiliary
Euclidean structure on $T_x M$.

\begin{definition}
A {\it Finsler metric\/} on a manifold $M$ is a function
$L : TM \rightarrow [0,\infty)$ that is
\begin{itemize}
\item positive and smooth outside the zero section;
\item positively homogeneous of degree one;
\item for each $x \in M$ the unit tangent sphere in $T_xM$
      is quadratically convex.
\end{itemize}
The metric $L$ is said to be {\it reversible\/} if for any tangent vector $v$,
we have that $L(v) = L(-v)$.
\end{definition}

Among the examples of reversible Finsler metrics we find Riemannian metrics,
submanifolds of Minkowski spaces (i.e., normed spaces whose unit spheres are
quadratically convex), and the  Hilbert geometries. Some examples of non
reversible Finsler metrics are the Katok examples
(see \cite{Katok} and \cite{Ziller}) of Finsler spheres with only two closed
geodesics, Bryant's examples of Finsler metrics on the 2-sphere
with constant curvature (\cite{Bryant:1, Bryant:2}), and the image of any
Riemannian metric by a small homogeneous canonical transformation.

In many cases, we will prefer to work with the duals of Finsler metrics:

\begin{definition}
An {\it optical Hamiltonian\/} on (the cotangent bundle of) a manifold $M$ is a
function $H : T^{*}M \rightarrow [0,\infty)$
that is
\begin{itemize}
\item positive and smooth outside the zero section;
\item positively homogeneous of degree one;
\item for every $x \in M$ the unit cotangent sphere in $T_{x}^* M$
      is quadratically convex.
\end{itemize}
The Hamiltonian $H$ is said to be {\it reversible\/} if for any covector $p$,
we have that $H(p) = H(-p)$.
\end{definition}

The duality between Finsler metrics and optical Hamiltonians is given by
the Legendre transform. This geometric transformation is best explained
on a single vector space:

Let $V$ be a finite-dimensional real vector space and let $S \subset V$ be a
quadratically convex hypersurface enclosing the origin. If $v$ is a point in $S$,
there is a unique covector $\xi \in V^*$ such that the hyperplane
$\xi = 1$ is tangent to $S$ at $v$ and the half-space $\xi \leq 1$ contains
$S$. The map that sends $v$ to $\xi$ is called the {\it Legendre transform\/}
and will be denoted by
$$
{\mathcal L} : S \longrightarrow V^* \ .
$$
The image of $S$ under $\mathcal L$, the {\it dual\/} of $S$, is denoted
by $S^*$ and is also a quadratically convex hypersurface that encloses
the origin. If we now use $S^* \subset V^*$ to define the  Legendre transform
$$
{\mathcal L}^* : S^* \longrightarrow V \  ,
$$
it is easy to see that $(S^*)^* = S$ and that
$$
{\mathcal L}({\mathcal L}^* (\xi)) = \xi
\quad \mbox{and}  \quad
{\mathcal L}^*({\mathcal L}(v)) = v \ .
$$

If $L$ is a Finsler metric on a manifold $M$, we can perform the above
construction on each tangent space $T_x M$, $x \in M$, and define a
diffeomorphism between $TM \setminus 0$ and $\ctg{M} \setminus 0$ which
we shall still call the Legendre transform and denote by $\mathcal L$.
With this notation, the function $H := L \circ {\mathcal L}^{-1}$
is an optical Hamiltonian. Conversely, if $H$ is an optical Hamiltonian
and ${\mathcal L}^*$ is its Legendre transform,
$L := H \circ {\mathcal L}^{*\, -1}$ is a Finsler metric. Using the standard
terminology from classical mechanics, we shall say that $L$ is the Lagrangian
of $H$, and $H$ is the Hamiltonian of $L$.

Given a Finsler metric $L$ on a manifold $M$, we define the length of a smooth
curve $\gamma : [a,b] \rightarrow M$ by
$$
\mbox{length}(\gamma) := \int_{a}^{b} L(\dot{\gamma}(t)) dt \ .
$$
When $L$ is reversible, this defines a length structure and a metric on $M$.
In general the length of a curve is only invariant under reparameterizations
that preserve the orientation.

The {\it geodesics\/} of a Finsler manifold $(M,L)$ are those curves that satisfy
the Euler-Lagrange equations for the Lagrangian $L$. In the reversible case
geodesics locally minimize length, and in the non reversible case they locally
minimize oriented length.

If $H$ is an optical Hamiltonian, we denote the sublevel set
$H \leq 1$ by $\cod{(M,H)}$ and its boundary, the level surface
$H = 1$, by $\csph(M,H)$. On this surface we use the canonical 1-form
$\alpha$ to define the {\it Reeb vector field\/} $X_{H}$ by the equations
$$
d\alpha(X_{H},\cdot) = 0 \quad \mbox{and} \quad
\alpha(X_{H}) = 1 \ .
$$

The integral curves of the Reeb vector field are usually called the
{\it characteristics\/} of $\csph{(M,H)}$.
A basic fact that we will use in the next two sections is that if
$H$ is an optical Hamiltonian and  $\gamma$ is a characteristic, then the
projection of $\gamma$ to $M$ is a geodesic for the associated Finsler metric,
$L$. Moreover, the image of $\gamma$ under the Legendre transform of $H$ is
the velocity curve of this  geodesic. Conversely if $c$ is a geodesic, then
the image of $\dot{c}$ under the Legendre transform of $L$ is a characteristic
of $\csph{(M,H)}$. The length of $c$ can be computed ``upstairs" as
the {\it action\/} of $\gamma$:
$$
\mbox{action}(\gamma) := \int_{\gamma} \alpha \ .
$$

Notice that if a homogeneous canonical transformation
$\phi : \ctg{M} \setminus 0 \rightarrow \ctg{M} \setminus 0$ preserves
$\csph{(M,H)}$, then it sends characteristics to characteristics and preserves
their action.

One of the advantages of the Hamiltonian viewpoint in Finsler geometry, is that
it suggests a natural definition of volume:

\begin{definition}
Let $L$ be a Finsler metric on an $n$-dimensional manifold $M$ and let
$H$ be its Hamiltonian. The {\it Holmes-Thompson volume\/} of $(M,L)$ is
defined as the symplectic volume of the set
$$
D^{*}(M,H) := \{p \in \ctg{M} : H(p) \leq 1\}
$$
divided by the volume of the Euclidean unit ball of dimension $n$. The
$k$-{\it area\/} of a $k$-dimensional submanifold of $M$ is defined as the
Holmes-Thompson volume of the submanifold with its induced Finsler metric.
\end{definition}

This definition, with its connections to symplectic geometry, convex geometry
(\cite{Holmes-Thompson,Thompson}), the Fourier transform
(\cite{Weil,Alvarez-Fernandes:2}), and integral geometry
(\cite{Schneider-Wieacker,Alvarez-Fernandes:1,Schneider}), has marked
advantages over the Hausdorff measure. Nevertheless, there is a simple
relationship between these two important notions of volume:

\begin{theorem}[Dur\'an, \cite{Duran}]\label{duran}
Let $(M,L)$ be a Finsler manifold with finite Hausdorff measure. The
Holmes-Thompson volume of $(M,L)$ is less than or equal to its Hausdorff
measure. Equality holds if and only if the metric is Riemannian.
\end{theorem}

Note that an optical Hamiltonian on a compact manifold is a star Hamiltonian,
and that $D^{*}(M,H)$ is a star body. This allows us to use the theory of dual
mixed volumes to define relative invariants of Finsler metrics:

\begin{definition}
Let $L_{1}$ and $L_{2}$ be two Finsler metrics on a compact manifold $M$ and let
$H_{1}$ and $H_{2}$ be their respective Hamiltonians. We define the $k$-th dual
mixed volume of $L_{1}$ and $L_{2}$ as the quantity
$$
\dmv_{k}(M;L_{1},L_{2}) := \dmv_{k}(D^{*}(M,H_{1}),D^{*}(M,H_{2})) \ .
$$
\end{definition}

Whenever we consider the first metric as a model metric against which
other metrics are to be compared (for example, an invariant metric on a
homogeneous space) we shall denote it by $L_{0}$ and define
$$
\emv_{k}(M,L) := \frac{1}{V(D^{*}(M,H_{0}))} \dmv_{k}(M;L,L_{0}) \ .
$$

Note that the main inequality and the definition of the Holmes-Thompson volume
easily imply the following result:

\begin{proposition}\label{basic}
If $L$ and $L_{0}$ are Finsler metrics on an $n$-dimensional compact manifold
$M$, then
$$
\emv_{k}(M,L)^{n} \leq \frac{\vol(M,L)^{n-k}}{\vol(M,L_{0})^{n-k}} \ .
$$
\end{proposition}

\section{Finslerian extensions of Pu's theorem}\label{inequalities}

In this section, we apply the theory of dual mixed volumes to extend the
works of Loewner and Pu on isosystolic inequalities to the Finsler setting.

The $k$-th {\it systole\/} of a Finsler manifold $(M,L)$,
denoted by $\sys_{k}(M,L)$, is defined the infimum of the  volumes of all
$k$-dimensional submanifolds not homologous to zero. When $k = 1$, it is usual
to change the definition to the infimum of the lengths of all
non-contractible curves on $M$.

For star Hamiltonians on cotangent bundles of manifolds that are not simply
connected, we define the 1-systole as the infimum of the actions of all
closed characteristics of $\csph(M,H)$ that project to non-contractible curves
on $M$. A theorem of Cieliebak (Theorem~2 in Part II of \cite{Cieliebak})
guarantees the existence of such characteristics for smooth star Hamiltonians on a
multiply-connected compact manifold.

As a rule, sharp isosystolic inequalities are only known for metrics in certain
conformal classes. The methods in the next two sections allow us to go somewhat
further, but in the present section all the results will be connected with
notions of conformality. The first of these notions is a straight-forward
generalization of conformality for Riemannian metrics.

\begin{definition}
Two Finsler metrics $L_{1}$ and $L_{2}$ on a manifold $M$ are said to be
{\it conformal\/} if there exists a diffeomorphism $\varphi : M \rightarrow M$
and a smooth positive function $\rho$ on $M$ such that
$\rho L_{2} =  \varphi^{*}L_{1}$.
\end{definition}

Since $1$-systoles and volumes are invariant under homogeneous
canonical transformations that are isotopic to the identity, any sharp
isosystolic inequality we prove for a class of Finsler metrics  will
hold for all star Hamiltonians obtained from the Hamiltonians of these
metrics by composition with some homogeneous canonical transformation that
is isotopic to the identity.  A somewhat more subtle notion of conformality
that involves homogeneous canonical transformations is as follows:

\begin{definition}
Two star Hamiltonians $H_{1}$ and $H_{2}$ on the cotangent of a manifold $M$
are said to be {\it s-conformal\/} if there exists a homogeneous canonical
transformation $\phi : \ctg{M} \setminus 0 \rightarrow  \ctg{M} \setminus 0$
and a smooth positive function $\rho$ on $M$
such that  $\rho H_{2} = H_{1} \circ \phi$.
\end{definition}

Note the slight abuse of notation in identifying the function $\rho$ with
its pull-back to $\ctg{M}$.

The definition of s-conformal Hamiltonians is not useful in Riemannian
geometry since even small homogeneous canonical transformations do not
generally send Riemannian metrics to Riemannian metrics. However, in the
Finsler and Hamiltonian setting, this notion allows us to recognize that
hidden symmetries can play a role in the proof of sharp isosystolic
inequalities.

The two main results of this section are:

\begin{theorem}\label{inequality:1}
Let $M$ be a compact homogeneous space and let $L_{0}$ be an invariant
Finsler metric on $M$. If $L$ is a Finsler metric conformal to $L_{0}$,
then
$$
\frac{\sys_{k}^{n}(M,L)}{\vol^{k}(M,L)} \leq
\frac{\sys_{k}^{n}(M,L_{0})}{\vol^{k}(M,L_{0})} .
$$
Moreover, equality holds if and only if $L$ is isometric to a multiple of
$L_{0}$.
\end{theorem}

\begin{theorem}\label{inequality:2}
Let $L_{0}$ be the standard Riemannian metric of curvature one on
$\RP^n$ and let $H_{0}$ be its Hamiltonian. If $L$ is a Finsler metric on
$\RP^n$ whose Hamiltonian $H$ is s-conformal to $H_{0}$, then
$$
\frac{\sys_{1}^{n}(\RP^n,L)}{\vol(\RP^n,L)} \leq
\frac{\sys_{1}^{n}(\RP^n,L_{0})}{\vol(\RP^n,L_{0})} \ .
$$
Moreover, equality holds if and only if $H = H_0 \circ \phi$ for some
homogeneous canonical transformation $\phi$.
\end{theorem}

The idea of the proof theorem~\ref{inequality:1} is very simple: by
proposition~\ref{basic}, we know
that
$$
\emv_{n-k}(M,L)^{n} \leq \frac{\vol(M,L)^{k}}{\vol(M,L_{0})^{k}} \ .
$$
Since $k$-systoles  and volumes are invariant under isometries,
we may assume that $L = \rho L_{0}$ with $\rho$ a smooth
positive function on $M$. The proof reduces to proving the inequality
$$
\emv_{n-k}(M,\rho L_{0}) \geq
\frac{\sys_{k}(M, \rho L_{0})}{\sys_{k}(M,L_{0})} \ .
$$
This hinges on interpreting $\emv_{n-k}(M,L)$ as different averages of the
radial function. For this we shall need two trivial lemmas:

\begin{lemma}\label{average:one}
Let $L_{0}$ be a Finsler metric on an $n$-dimensional manifold $M$.
If $\rho$ is a smooth positive on $M$, the quantity $\emv_{n-k}(M,\rho L_{0})$
is the average of the $k$-th power of $\rho$ over the manifold $M$:
$$
\emv_{n-k}(M,\rho L_{0}) = \frac{1}{\vol(M,L_{0})} \int_{M} \rho^{k} \ dV^{0} \ ,
$$
where $dV^{0}$ is the density for the Holmes-Thompson volume on $(M,L_{0})$.
\end{lemma}

\begin{lemma}\label{average:two}
Let $G$ be a compact Lie group and let $\mu$ be the Haar measure on $G$ normalized
so that the measure of $G$ equals one. If $Q$ is a compact manifold on which
$G$ acts transitively and $\nu$ is an invariant measure on $Q$, then for any
function  $f \in L^{1}(Q,\nu)$, we have that
$$
\int_{G} f(g\cdot x) \ d\mu = \frac{1}{\nu(Q)} \int_{Q} f(y) \ d\nu \ .
$$
In other words, the average of the pullback of $f$ to $G$ equals the average
of $f$ on $Q$.
\end{lemma}

Putting both lemmas together, we have the following proposition:

\begin{proposition}\label{average:three}
Let $M$ be homogeneous space under the left-action of a compact Lie
group $G$ and let $\mu$ be the Haar measure on $G$ normalized so that the
measure of $G$ equals one. If $L_{0}$ is an invariant Finsler metric on $M$
and $\rho$ is a smooth positive function, then
$$
\emv_{n-k}(M,\rho L_{0}) = \int_{G} \rho^{k}(g \cdot x) \ d\mu \
$$
for any $x \in M$.
\end{proposition}

\proof
Applying lemma~\ref{average:two} with $Q := M$, $f := \rho^{k}$, and
$\nu := dV^{0}$, the volume density of the Holmes-Thompson volume of
$(M,L_{0})$, we obtain that
$$
\int_{G} \rho^{k}(g \cdot x) \ d\mu =
\frac{1}{\vol(M,L_{0})} \int_{M} \rho^{k} \ dV^{0} .
$$
By lemma~\ref{average:one}, this quantity equals $\emv_{n-k}(M,\rho L_{0})$.
\qed

\begin{proposition}\label{conformal}
Let $M$ be a compact homogeneous space and let the model metric $L_{0}$ be
an invariant Finsler metric on $M$. If $\rho$ is a smooth positive function
on $M$, then
$$
\emv_{n-k}(M,\rho L_{0}) \geq \frac{\sys_{k}(M,\rho L_{0})}{\sys_{k}(M,L_{0})} \ .
$$
\end{proposition}

\proof
By proposition \ref{average:three}, we have that
$$
\emv_{n-k}(M,\rho L_{0}) = \int_{G} \rho^{k}(g \cdot x) \ d\mu \
$$
for any point $x \in M$.

If $N_{i} \subset M$ is a sequence of $k$-dimensional submanifolds
not homologous to zero whose $k$-dimensional volumes decrease to
$\sys_{k}(M,\rho L_{0})$,
 \begin{eqnarray*}
\emv_{n-k}(M,\rho L_{0}) \ \vol_{k}(N_{i},L_{0}) & = &
\int_{N} \left(\int_{G} \rho^{k} \ d\mu \right) \ dV^{0}_{k} \\
& = & \int_{G} \left( \int_{g(N_{i})} \rho^{k} \ dV^{0}_{k} \right) \ d\mu .
\end{eqnarray*}
Since $\rho^{k} \ dV^{0}_{k}$ is the $k$-area density of the Finsler metric
$\rho L_{0}$ and $g(N_{i})$ is not homologous to zero, the last integral is
greater than or equal than the $k$-systole of $(M,\rho L_{0})$. Therefore, for all
$i$ we have that
$$
\emv_{n-k}(M,\rho L_{0}) \ \vol_{k}(N_{i},L_{0}) \geq \sys_{k}(M,\rho L_{0}) .
$$
Taking the limit as $i$ tends to infinity yields the desired inequality.
\qed

The proof of theorem~\ref{inequality:2} is somewhat more subtle because we
need to average with respect to hidden symmetries. Like in the case of
theorem~\ref{inequality:1}, everything boils down to proving the following
result:

\begin{proposition}
Let $L_{0}$ be the standard Riemannian metric of curvature one on
$\RP^n$ and let $H_{0}$ be its Hamiltonian. If $L$ is a Finsler metric on
$\RP^n$ whose Hamiltonian $H$ is s-conformal to $H_{0}$, then
$$
\emv_{n-1}(M,L) \geq \frac{\sys_{1}(M,L)}{\sys_{1}(M,L_{0})} \ .
$$
\end{proposition}

\proof
Since $H$ is s-conformal to $H_{0}$, there exists a smooth positive function
$\rho$ on $M$ and a homogeneous canonical transformation $\phi$ such that
$\rho H = H_{0} \circ \phi$. Notice the slight abuse notation in denoting
$\rho$ and its pullback to the cotangent bundle by the same symbol.

If we lift the standard action of $SO(n+1)$ on
$\RP^n$ to $T^{*}\RP^n$ and conjugate by $\phi$ we obtain a left-action of
$$
SO(n+1) \times (\ctg{\RP^n} \setminus 0) \longrightarrow (\ctg{\RP^n} \setminus 0)
$$
by homogeneous canonical transformations that preserve the Hamiltonian
$H_{0} \circ \phi =: K_{0}$. Notice that the geodesic flow of $K_{0}$ is
symplectically conjugate to that of $H_{0}$ and, therefore, all geodesics are
closed of length $\pi$.

By lemma~\ref{average:two}, we have that
$$
\emv_{n-1}(\RP^n,L) = \int_{SO(n+1)} \rho (g \cdot p) \ d\mu \ .
$$
Mimicking the proof of proposition~\ref{conformal}, we let $\gamma$
be a closed characteristic of $\csph{(\RP^n,K_{0})}$ and write
\begin{eqnarray*}
\emv_{n-1}(\RP^n,L) \ \sys_{1}(\RP^n,L_{0}) & = &
\int_{\gamma} \left(\int_{SO(n+1)} \rho \ d\mu \right) \, \alpha \\
& = & \int_{SO(n+1)} \left( \int_{g(\gamma)} \rho \, \alpha \right) \ d\mu .
\end{eqnarray*}

We would like to argue that the integral of $\rho \, \alpha$ along $g(\gamma)$
is the length of some non-contractible curve in the metric $L$, and
conclude that it must be greater than $\sys_{1}(M,L)$. In order to do this
we must show that the curve $g(\gamma)$ is the image under the Legendre
transform of $L$ of the velocity curve of some non-contractible
curve in $\RP^n$.  Since, as a rule, homogeneous canonical transformations mix
momentum and position, this is not entirely obvious.

We start by noticing that since $g$ preserves both $K_{0}$ and the
canonical 1-form $\alpha$, $g(\gamma)$ is also a characteristic of
$\csph{(\RP^n,K_{0})}$, and, hence, its image under the Legendre transform
of $K_{0}$ is the velocity curve of some curve in $c$ in $\RP^n$. Since $SO(n+1)$
is connected, the curve $c$ cannot be contractible. Now, we use that
$\rho H = K_{0}$ and that $\rho$ is constant on the fibers to conclude that
the Legendre transform of $H$ sends the curve $g(\gamma)$ to the velocity
curve of a suitable reparameterization of $c$.

We conclude that
\begin{eqnarray*}
\emv_{n-1}(\RP^n,L) \ \sys_{1}(\RP^n,L_{0}) & = &
\int_{SO(n+1)} \left( \int_{g(\gamma)} \rho^{k} \ \alpha \right) \ d\mu \\
& \geq & \sys_{1}(\RP^n,L) \ .
\end{eqnarray*}
\qed

In the Riemannian case, isosystolic inequalities on the two-dimensional torus
and the projective plane are substantially strengthened by the fact that
any Riemannian metric on these surfaces is conformal to a homogeneous metric.
This is utterly false in the Finsler case, nevertheless it seems that
any {\it reversible\/} Finsler metric on $\RP^2$ is s-conformal to the
standard Riemannian metric of curvature one:

\begin{conjecture}
If $L$ is a reversible Finsler metric on the projective plane, $\RP^2$,
then there exists a smooth positive function of $\RP^2$ such that
the Finsler metric $\rho L$ has periodic geodesic flow.
\end{conjecture}

Using the results of this section, the uniformization conjecture would
lead to an alternate proof of Ivanov's generalization of Pu's theorem.

\begin{theorem}[Ivanov, \cite{Ivanov}]\label{Finsler-Pu}
Any reversible Finsler metric $L$ on the projective plane satisfies
the inequality
$$
\frac{2}{\pi}\sys_{1}^{2}(\RP^2,L) \leq \vol(\RP^2,L) \ .
$$
Equality holds if and only if the geodesic flow of $L$ is periodic.
\end{theorem}

\section{Isosystolic inequalities for commuting Hamiltonians}

In this section, we move away from notions of conformality and establish
isosystolic inequalities under a different set of hypotheses. Namely, we
shall consider star Hamiltonians with periodic flow and
star Hamiltonians commuting with them.

\begin{theorem}\label{periodic}
If $H_{0}$ is a smooth star Hamiltonian on $\ctg{\RP^n}$ with periodic flow,
all simple characteristics on $H_0 = 1$ project to closed
non-contractible curves in $\RP^n$ and
$$
\frac{\sys_{1}^{n}(\RP^n,H_{0})}{\vol(\RP^n,H_{0})}
= \frac{2\pi^n}{(n+1)\varepsilon_{n+1}} \ ,
$$
where $\varepsilon_{k}$ is the volume of the Euclidean unit ball of dimension
$k$.
\end{theorem}

\proof
The periodicity of the flow implies that all simple characteristics
are homotopic. Since by  Theorem~2, Part II, of \cite{Cieliebak}, some closed
characteristic on $H_0 = 1$ projects to a non-contractible curve in
$\RP^n$, all characteristics project to non-contractible curves.

The proof of equality relies on the results and techniques
of Weinstein (\cite{Weinstein}) and Yang (\cite{Yang}). They proved that
the volume of a Riemannian metric on the $n$-sphere whose
geodesic flow is periodic with period $\ell$ is equal to the volume of the
Euclidean $n$-sphere of radius $\ell/2\pi$. Their techniques are symplectic and
topological, and carry over without modification to smooth star Hamiltonians
on the cotangent bundle of the sphere or real projective space.
\qed

\begin{theorem}\label{inequality:3}
If $H_{0}$ is a smooth star Hamiltonian on $\ctg{\RP^n}$ with periodic flow and
$H$ is a star Hamiltonian that is constant along the orbits of $X_{H_{0}}$,
then
$$
\frac{\sys_{1}^{n}(\RP^n,H)}{\vol(\RP^n,H)} \leq
\frac{\sys_{1}^{n}(\RP^n,H_{0})}{\vol(\RP^n,H_{0})} \ .
$$
Moreover, equality holds if and only if $H$ is a fixed multiple of
$H_{0}$.
\end{theorem}

Since for any star body $A \subset T^* M$
$$
\emv_{n-1}(A) = \frac{1}{V(A)} \int_{\partial U} \rho_A \Omega
\geq \min \rho_A \ ,
$$
the proof of the theorem above reduces to that of the following inequality:

\begin{proposition}\label{minumum}
Let $H$ and $H_{0}$ be as in theorem~\ref{inequality:3}. If $\rho$ denotes the
radial function of the star body $H \leq 1$, then
$$
\min \rho \geq \frac{\sys_1(\RP^n,H)}{\sys_1(\RP^n, H_{0})} \ .
$$
\end{proposition}

The gist of the proof is to show that if  $\lambda := \min \rho$ is
attained at a point $p_{0} \in \csph{(M,H_{0})}$, then the integral curve of
$X_{H}$ on $\csph{(M,H)}$ with initial condition $p := \rho(p_{0}) \, p_{0}$ is
closed, projects to a non-contractible curve, and its action equals
$\lambda \, \sys_1(\RP^n,H_{0})$.

\begin{lemma}
Consider the map $T : \csph{(\RP^n,H_{0})} \rightarrow \csph{(\RP^n,H)}$ defined
by $p \mapsto \rho(p) \, p$.  If $p$ is a critical point of $\rho$,
then $T_{*}(X_{H_{0}}(p))$ is a positive multiple of $X_{H}(\rho(p) \, p)$
\end{lemma}

\proof
Since
$$
T^*\left(d\alpha \lfloor T_{*}(X_{H_{0}})\right) =
(d\rho \wedge \alpha) \lfloor X_{H_{0}} \ ,
$$
we have that when $d\rho(p) = 0$, $d\alpha \lfloor T_{*}(X_{H_{0}})$ must be
zero. The kernel of $d\alpha$ on $\csph(\RP^n,H)$ is one dimensional and, hence,
$T_{*}(X_{H_{0}}(p))$ is a multiple of $X_{H}(\rho(p) \, p)$. To see that it
is a positive multiple we remark that
$$
\alpha(T_{*}(X_{H_{0}})) = (T^*\alpha)(X_{H_{0}}) = \rho \alpha(X_{H_{0}}) =
\rho \ ,
$$
which is always positive.
\qed

\noindent
{\it Proof of proposition~\ref{minumum}.\ }
We now proceed to build a closed characteristic on $\csph{(\RP^n,H)}$ that
is closed, projects to a non-contractible curve in $\RP^n$, and such that
its action equals $\lambda \, \sys(\RP^n,H_{0})$.

Let $p_{0}$ be a point in $\csph{(\RP^n,H_{0})}$ where the minimum of $\rho$
is attained, and let $\sigma(t)$ be a closed characteristic with
$\sigma(0) = p_{0}$. Notice that since $\rho$ is constant along characteristics,
$\rho(\sigma(t))$ is constant and all the points $\sigma(t)$ are critical points
of $\rho$. This, together with the fact that $T_{*}(X_{H_{0}}(p))$ is a positive
multiple of $X_{H}(\rho(p) \, p)$ whenever $p$ is a critical point, implies that
$$
\gamma(t) := T(\sigma(t)) = \lambda \sigma(t)
$$
is a characteristic of $\csph{(\RP^n,H)}$. Moreover, $\gamma$ projects down to
the same curve as $\sigma$. By Theorem~\ref{periodic}, this implies that
$\gamma$ projects to to a non-contractible curve on $\RP^n$.

It is easy to compute the action of $\gamma$:
\begin{eqnarray*}
\int_{\gamma} \alpha & = & \int_{T(\sigma)} \alpha = \int_{\sigma} T^* \alpha \\
                     & = & \int_{\sigma} \lambda \alpha =
                           \lambda \, \sys(\RP^n,H_{0}) \ .
\end{eqnarray*}
\qed

Identifying $SO(3)$ with $\RP^3$ and recalling (see \cite{Abraham-Marsden})
that left-invariant Hamiltonians poisson-commute with the Hamiltonian of
the bi-invariant metric on $SO(3)$, we have the following corollary:

\begin{corollary}
If $L$ a Finsler metric on $SO(3)$ that is conformal to
a left-invariant Finsler metric, then
$$
\frac{\sys_{1}(SO(3),L)^3}{\vol(SO(3),L)} \leq \pi \ .
$$
Equality holds if and only if $L$ is bi-invariant.
\end{corollary}

\section{Hamiltonian averaging and isosystolic inequalities}\label{averaging}

The main result of this section is a Hamiltonian generalization of Berger's
infinitesimal isosystolic inequality for Riemannian metrics on real projective
spaces (\cite{Berger}).

\begin{theorem}\label{inequality:4}
If $H_{t}$ is a smooth path of smooth star Hamiltonians on $\ctg{\RP^n}$ such
that the flow of $H_{0}$ is periodic, then there exists another
smooth path of smooth star Hamiltonians, $K_t$, that agrees to first order
with $H_t$ at $t = 0$ and satisfies the isosystolic inequality
$$
\frac{\sys_{1}^{n}(\RP^n,K_t)}{\vol(\RP^n,K_t)} \leq
\frac{\sys_{1}^{n}(\RP^n,H_{0})}{\vol(\RP^n,H_{0})} \ .
$$
\end{theorem}

For the  proof we will need the following mild adaptation of a  classical
(and easy) result in the theory of normal forms of Hamiltonian systems
(see lemma 3.3 in \cite{Cushman}):

\begin{lemma}\label{normal_form}
Let $H_{0}$ be a star Hamiltonian on $\ctg{\RP^n}$ whose flow is periodic.
Any Hamiltonian $H$ on $T^*\RP^n \setminus 0$ that is homogeneous of degree
one can be written in a unique way as $E + \{H_{0},F\}$, where $E$ and $F$
are also homogeneous Hamiltonians of degree one and $\{H_{0},E\} = 0$.
\end{lemma}

\noindent
{\it Proof of theorem~\ref{inequality:4}.\ }
By theorem~\ref{inequality:3}, we know that all star Hamiltonians  of
the form $H_{0} + tE$ with $\{H_{0},E\} = 0$ satisfy the isosystolic inequality.
Moreover, If $F$ is a Hamiltonian that is homogeneous of degree one, and
$\phi_{t}$ is its flow, then for sufficiently small $t$, the function
$(H_{0} + tE) \circ \phi_{t} =: K_{t}$ is a star Hamiltonian and
satisfies the isosystolic inequality.

Writing $H_{t} = H_{0} + tH_{1} + O(t^2)$, we see that it is
possible to find $E$ and $F$ such that $K_{t}$ and $H_{t}$ agree up to order
one at $t = 0$. Indeed
$$
\frac{d K_{t}}{dt}(0) = E + \{H_{0},F\} \ ,
$$
and, according to lemma~\ref{normal_form}, we may always find $E$ and $F$
such that $E + \{H_{0},F\} = H_{1}$.
\qed


\end{document}